\documentclass{article}
\usepackage{amsmath}
\usepackage{amssymb}
\usepackage{amsthm}
\usepackage{tikz}

\begin{document}

\title{Backward bifurcations and multistationarity}

\author{Alexis Nangue\\
Higher Teachers' Training College, University of Maroua\\
P.O. Box 55, Maroua, Cameroon\\  
and\\
Alan D. Rendall\\
Institut f\"ur Mathematik\\
Johannes Gutenberg-Universit\"at, Staudingerweg 9\\
D-55099 Mainz, Germany}

\date{}

\maketitle

\begin{abstract}
The theory of backward bifurcations provides a criterion for the existence of
positive steady states in epidemiological models with parameters where the
basic reproductive ratio is less than one. It is often seen in simulations
that this phenomenon is accompanied by multistationarity, i.e. the existence of
more than one positive steady state, but the latter circumstance is not implied
by the general theory. The central result of this paper is a theorem which gives
a criterion for the existence of one stable and one unstable positive steady
state for parameters where the basic reproductive ratio is less than one. It
also gives a criterion for the existence of one stable and one unstable positive
steady state in the case that the basic reproductive ratio is greater than one.
These steady states arise in a bifurcation. It is shown that in one case, a
model for the in-host dynamics of hepatitis C, this result can be used as a
basis for showing the existence of parameters for which there are two positive
steady states, one of which is stable. Thus, in particular, multistationarity
is proved in that case. It is also shown to what extent this theorem can be
applied to some other models which have been studied in the literature and
what new results can be obtained. In that context the new approach is compared
with those previously known.
\end{abstract}  

\section{Introduction}

In epidemiology a central role is played by the basic reproductive ratio
${\cal R}_0$ \cite{diekmann00}. In the simplest models an endemic equilibrium,
i.e. a positive steady state, exists for ${\cal R}_0>1$ and does not exist for
${\cal R}_0<1$. There are, however, models with more complicated behaviour
where a positive steady state does exist in some situations with
${\cal R}_0<1$. This phenomenon has been referred to as a backward
bifurcation. More information concerning this concept can be found in
\cite{martcheva19}. A common visual representation of a backward bifurcation
is given by plotting the value of a variable $I$ corresponding to the size of
an infected population at steady state against ${\cal R}_0$. The result is a
curve which passes through the point $({\cal R}_0,I)=(1,0)$. In the case of a
backward bifurcation its slope at that point is negative and close to that
point the positive steady state is unstable. This contrasts with the more
familiar case of a forward bifurcation where the slope at the point $(1,0)$ of
the curve defined by steady states is positive and close to that point the
steady state is stable. A general sufficient condition for the occurrence of a
backward bifurcation and a proof that under this condition the properties
mentioned above hold have been given in Theorem 4 of \cite{vandendriessche02}.

Plots of the curve defined by steady states in a backward bifurcation obtained
by simulation often show features additional to those discussed up to now. It
is often seen that there is a point where the curve turns around, with its
slope becoming positive. Up to the turning point the solution is unstable and
after that point it is stable. There is a fold (saddle-node) bifurcation
\cite{perko01} at the turning point. To our knowledge there is no analytical
proof in the literature of a general criterion for the occurrence of these
additional phenomena. One special case where more has been proved is that
where the steady states are in one-to-one correspondence with the roots of a
quadratic equation whose coefficients are functions of the parameters
occurring in the system. For this see \cite{brauer04} and section 2.2 of
\cite{martcheva19}. The aim of what follows is to provide a general criterion
for the occurrence of a fold bifurcation of this type and the associated
existence of more than one positive steady state. This situation corresponds to
the occurrence of a certain type of bifurcation defined in the next section.
As will be explained later it is an unfolding of a pitchfork bifurcation.
This criterion is then applied to some examples. In fact this type of
bifurcation can give rise to the appearance of two steady states either in the
case ${\cal R}_0<1$ just discussed or in the case ${\cal R}_0>1$. A criterion
is given for distinguishing between these two cases.

The paper is organized as follows. Theorem 1 of Section \ref{main} is the
central result of the paper. It provides general criteria for the existence
of parameters in a system of ODE arising in epidemiology for which
there exist a certain number of steady states with certain stability
properties. Some of these steady states bifurcate from a steady state on
the boundary of the positive orthant so that their properties can be explored
by a local analysis. In Section \ref{brauer} the results of this paper are
compared with similar results on backward bifurcations for a certain
vaccination model which have been obtained by more elementary methods in
\cite{brauer04}. In Section \ref{martcheva} it is shown that results obtained
in \cite{martcheva19} on the existence of more than one steady state in a model
for cholera with vaccination introduced there can be strengthened by an
application of Theorem 1 to give a rigorous proof of their stability
properties. Section \ref{hepc} provides an analysis of an in-host model of
hepatitis C previously studied by the authors in \cite{nangue22} and
\cite{nangue23}. By applying Theorem 1 we obtain a new result for this
model. It says that there are parameters with ${\cal R}_0<1$ for which there
exist exactly two positive steady states, one stable and one unstable. These
arise in a backward bifurcation. In particular this proves that it is possible
to have a stable steady state in the case ${\cal R}_0<1$. For this system there
exist at most three steady states for any choice of the parameters. It is also
shown that this system has a limiting case where there are parameters for which
there exists a whole continuum of steady states. The last section gives
conclusions and an outlook on possible future extensions of this work.

\section{The main theorem}\label{main}

The formulation of the theorem proved in this section is analogous to that
of Theorem 4 of \cite{vandendriessche02} and hence we start by recalling some
aspects of that paper. Theorem 4 of \cite{vandendriessche02} is formulated in
terms of a system of ODE $\dot x=f(x,\mu)$ depending on a parameter $\mu$. For
our theorem we require a second parameter and we consider a system of the form
$\dot x=f(x,\alpha_1,\alpha_2)$. Here $\alpha_1$ plays the role of the parameter
$\mu$ in \cite{vandendriessche02}. The unknown $x$ has $n$ components
$x_1,\ldots,x_n$. In \cite{vandendriessche02} a subset of these,
$x_1,\ldots,x_m$, is chosen as the set of populations of infected individuals.
We call these infected variables. A given system of ODE may have more than one
interpretation as an epidemiological model and there may be more than one
reasonable choice of the set of infected variables. In general the basic
reproductive ratio ${\cal R}_0$ depends on this choice but the sign of
${\cal R}_0-1$ does not. There is a point $x_0$ called the disease-free
equilibrium (DFE) on the boundary of the positive orthant which is a steady
state for all values of the parameter $\mu$. $X_s$ is the set of disease-free
states, defined by the equations $x_i=0$ for $i=1,2,\ldots,m$. The quantities
set to zero here are the infected variables and the DFE belongs to $X_s$. In
the formalism used in Theorem 4 of \cite{vandendriessche02} the variables are
ordered in such a way that the infected variables come first. When doing
calculations in examples it will in general be convenient to use a different
ordering of the variables. 

In \cite{vandendriessche02} the components $f_i$ of $f$ have a special
structure and are required to satisfy some conditions (A1)-(A5) defined in
terms of this structure. These details will not be repeated here and we will
only recall some key features of these conditions. $X_s$ is invariant as a
consequence of (A2) and (A4). It is assumed that $x_0\in X_s$ and that it is
stable within $X_s$. The right hand sides of the equations are of the form
$f_i={\cal F}_i-{\cal V}_i$. Here the parameter dependence has not been
written explicitly. The system $\dot x_i=-{\cal V}_i$ will be called the
truncated system. The set $X_s$ is invariant for the truncated system and
$x_0$ is a steady state for that system. (A5) is the condition that in the
truncated system the eigenvalues of the matrix $A=D_xf$ defining the
linearization about $x_0$ corresponding to eigenvectors lying in $X_s$ have
negative real parts. The matrix $A$ is of block
triangular form and the submatrix on the upper right is zero. Hence the set of
eigenvalues of the linearization is the union of the sets of eigenvalues of
the diagonal blocks $J_1$ and $J_4$. The condition (A5) implies that all
eigenvalues of $J_4$ have negative real parts. Hence the stability of $x_0$ is
determined by $J_1=F-V$ in the notation of \cite{vandendriessche02}. In that
paper it is proved that the sign of the maximum real part of an eigenvalue of
$J_1$ is equal to that of ${\cal R}_0-1$. For ${\cal R}_0=1$ the DFE has a
non-trivial centre manifold and in the general discussion of
\cite{vandendriessche02} it is assumed that that manifold is of dimension one.
It is also shown that the first $m$ components of the eigenvector corresponding
to the eigenvalue zero can be chosen to be non-negative. In
\cite{vandendriessche02} two quantities $a$ and $b$ are defined  which have
the property that the criterion for a backward bifurcation is that $a>0$ and
$b>0$. The criterion for a forward bifurcation is that $a<0$ and $b>0$. These
quantities will be discussed further below. 

To facilitate the description of the new theorem it is helpful to start by
giving an alternative account of the content of Theorem 4 of
\cite{vandendriessche02}. Consider a family of systems depending on a
parameter $\mu$ for which $x_0$ is always a steady state and ${\cal R}_0$
passes through one when $\mu$ passes through zero. Consider the extended system
where the equation $\dot\mu=0$ is adjoined to the original system. Assume that
${\cal R}_0$ is an increasing function of $\mu$ with ${\cal R}_0(0)=1$. We
call $(x_0,0)$ the bifurcation point. Consider the centre manifold of the
bifurcation point in the extended system. It is two-dimensional and the
restriction of the extended system to that manifold defines a one-dimensional
dynamical system depending on the parameter $\mu$. Let $u$ be a coordinate which
parametrizes the one-dimensional centre manifold of $x_0$, vanishes at $x_0$
and increases in the direction of the eigenvector tangent to the centre
manifold whose first $m$ components are non-negative. Because $u=0$ is a
steady state for all values of $\mu$ the restriction of the system to that
manifold is of the form $\dot u=f(u,\mu)$ with $f(u,\mu)=ug(u,\mu)$ for a
smooth function $g$.

In the situation considered in \cite{vandendriessche02} the bifurcation is a
non-degenerate transcritical bifurcation. Under the assumption that the steady
state is at zero for all values of the parameter this is characterized by the
conditions that $f_u(0,0)=0$, $f_{uu}(0,0)\ne 0$ and $f_{u\mu}(0,0)\ne 0$. This
is Theorem 3.12 of \cite{homburg20}. Note that the transcritical bifurcation
is also sometimes considered without the additional assumption that the steady
state is fixed. For this see Theorem 3.13 of \cite{homburg20} and the
discussion in \cite{perko01}. When $f$ can be expressed in terms of a function
$g$ as above the conditions on $f$ can be rewritten as conditions on $g$.
These are $g(0,0)=0$, $g_u(0,0)\ne 0$ and $g_\mu(0,0)\ne 0$. This means that
for $g$ the point $(0,0)$ is regular, i.e. not a bifurcation point. In the
situation usually considered in epidemiology the DFE is stable for
${\cal R}_0<1$ and unstable for ${\cal R}_0>1$. This means that when passing
through the bifurcation point in the direction of increasing ${\cal R}_0$,
which is also the direction of increasing $\mu$, the derivative $f_u$ must be
increasing. Hence it will be assumed from now on that $f_{u\mu}(0,0)>0$. The
parameters $a$ and $b$ defined in \cite{vandendriessche02} correspond to
$\frac12 f_{uu}(0,0)$ and $f_{u\mu}(0,0)$ respectively. To see this it is
important to understand the relation between conditions on the system on the
centre manifold and conditions on the full system from which it is derived by
restriction. The relevant facts are explained in \cite{kuznetsov10}, section
5.4. Under the given assumptions the linearization $A$ at the bifurcation
point has a one-dimensional kernel which is spanned by a vector $w$ and $A^T$
has a one-dimensional kernel which is spanned by a vector $v$. These vectors
can be normalized so that their inner product is one. Then
$a=\frac12\sum_{i,j,k=1}^nv_iw_jw_k
\frac{\partial ^2f_i}{\partial x_j\partial x_k}(x_0,0)$
and $b=\sum_{i,j=1}^nv_iw_j\frac{\partial ^2f_i}{\partial x_j\partial\mu}(x_0,0)$.
Note that while the quantity $a$ only depends on the system at $\mu=0$ the
quantity $b$ depends on the way the system varies with the parameter $\mu$.

In the theorem which follows the system $\dot x=f(x,\mu)$ is replaced by
$\dot x=f(x,\alpha_1,\alpha_2)$. Here $f$ is assumed to be smooth, i.e.
$C^\infty$. Setting $\mu=\alpha_1$ and regarding
$\alpha_2$ as fixed provides a definition of $a$ and $b$ for the latter
system. Adjoining the equations $\dot\alpha_1=0$ and $\dot\alpha_2=0$ to those
for $\dot x$ gives an extended system.

\noindent
{\bf Theorem 1} Consider the system $\dot x=f(x,\alpha_1,\alpha_2)$ for a
smooth function $f$ as an analogue of the disease transmission model of
\cite{vandendriessche02}. Suppose that for each fixed value of $\alpha_2$ the
system satisfies the conditions (A1)-(A5) of \cite{vandendriessche02} with
the same DFE for all values of $\alpha_1$. The equation on the centre
manifold of the extended system is of the form $\dot u=f(u,\alpha_1,\alpha_2)$.
Let $a=\frac12 f_{uu}(0,0,0)$ and $b= f_{u\alpha_1}(0,0,0)$ as in
\cite{vandendriessche02} and let $c=\frac13 f_{uuu}(0,0,0)$ and
$d=f_{uu\alpha_2}(0,0,0)$. Assume that for given parameter values $a=0$, $b>0$
and $c\ne 0$. Let $e$ be the value of $\frac{-bd+f_{uu\alpha_1}f_{u\alpha_2}}{2bc}$
at the bifurcation point. If $ce<0$ a backward bifurcation occurs for
$\alpha_2>0$ and for $ce>0$ a forward bifurcation occurs for $\alpha_2>0$. If
$c<0$ and $e>0$ there are parameters for which there are two positive steady
states with ${\cal R}_0<1$, one stable and one unstable. If $c>0$ and $e>0$
there are parameters for which there are two positive steady states with
${\cal R}_0>1$, one stable and one unstable.

\noindent
{\bf Proof} Note first that $f(u,\alpha_1,\alpha_2)=ug(u,\alpha_1,\alpha_2)$ for
a smooth function $g$. The conditions on $f$ translate to the conditions
$g_u(0,0,0)=0$, $g_{\alpha_1}(0,0,0)>0$, $g_{uu}(0,0,0)\ne 0$ and
$(-g_{\alpha_1}g_{u\alpha_2}+g_{\alpha_2}g_{u\alpha_1})(0,0,0)\ne 0$. This means in
particular that when $\alpha_2=0$
there is a generic fold bifurcation at the origin \cite{kuznetsov10}. Consider
the mapping $\Phi:(u,\alpha_1,\alpha_2)\mapsto (g(u,\alpha_1,\alpha_2),
g_u(u,\alpha_1,\alpha_2))$. Then the bifurcation points of the system with
$u\ne 0$ are the zeroes of $\Phi$. The aim is to apply the implicit function
theorem to solve the equation $\Phi=0$ for $u$ and $\alpha_1$ as functions of
$\alpha_2$. For this we have to check the invertibility of a $2\times 2$ matrix.
Since $g_u=0$ the determinant of this matrix is $-g_{uu}g_{\alpha_1}$ and this
quantity is indeed non-zero under the assumptions of the theorem. Thus there
exist unique smooth functions $U$ and $A_1$ in a neighbourhood of zero such
that $\Phi(U(\alpha_2),A_1(\alpha_2),\alpha_2)=0$. The derivatives of
these functions at zero are given by
\begin{eqnarray}
  &&  U'(0)=\left[-\frac{g_{u\alpha_2}}{g_{uu}}
     +\frac{g_{u\alpha_1}g_{\alpha_2}}{g_{uu}g_{\alpha_1}}\right](0,0,0),\\
&&A_1'(0)=
  -\frac{g_{\alpha_2}}{g_{\alpha_1}}(0,0,0).
\end{eqnarray}
(The quantity in brackets in the expression for $U'(0)$ is equal to $2e$.)
This means that close to the origin the positive bifurcation points lie on
a smooth curve. By continuity we have $g_{uu}\ne 0$ and $g_{\alpha_1}>0$ at those
points. Hence for each fixed value of $\alpha_2$ this is a generic fold
bifurcation and the zero set of $g$ can be approximated by a parabola. The
parabola opens to the right if $c<0$ and to the left if $c>0$. If
$e=U'(0)\ne 0$ then for $\alpha_2\ne 0$ we
have $u\ne 0$ at the bifurcation point. Hence the parabola intersects the
$\alpha_1$ axis transversally. If $ce<0$ the slope of the curve at its
intersection with the $\alpha_1$ axis is negative and we obtain
a backward bifurcation for $\alpha_2>0$ while if $ce>0$ the slope is positive
and we obtain a forward bifurcation for $\alpha_2>0$. The remaining statements
of the theorem follow by considering the relative positions of the parabola
and the line $u=0$. $\blacksquare$

When the conditions of this theorem are satisfied, i.e. $a=0$, $b>0$, $c\ne 0$
and $e\ne 0$ a fold bifurcation moves across the boundary of the positive
orthant as the parameter $\alpha_2$ is varied. This is analogous to the case
of a transcritical bifurcation, where a regular steady state moves across the
boundary of the positive orthant as a parameter is varied. Note that the sign
of $e$ can be reversed by reversing the sign of $\alpha_2$. The important
property of this quantity is that it be non-zero. Note that under the
hypotheses of the theorem a sufficient condition for this to be true is that
$f_{uu\alpha_1}(0,0,0)=0$ and $d=f_{uu\alpha_2}(0,0,0)\ne 0$. The conditions of
the theorem are expressed in terms of the restriction of the system to the
centre manifold and it is desirable to rewrite them in terms of the full
system. This can be done following the calculations of \cite{kuznetsov10}. The
result is
\begin{eqnarray}
  &&c=\frac13\sum_{i,j,k,l=1}^nv_iw_jw_kw_l\frac{\partial^3 f_i}{\partial x_j\partial x_k\partial x_l}(x_0,0,0)
      \nonumber\\
  &&-\sum_{i,j,k,l,p,q=1}^nv_iw_jw_kw_pB^q_l\frac{\partial^2 f_l}{\partial x_j\partial x_k}(x_0,0,0)
     \frac{\partial^2 f_i}{\partial x_q\partial x_p}(x_0,0,0),\\
&&d=\sum_{i,j,k=1}^nv_iw_jw_k\frac{\partial ^3f_i}
   {\partial x_j\partial x_k\partial\alpha_2}(x_0,0,0).
\end{eqnarray}
The expression for $c$, which corresponds to equation (8.128) of
\cite{kuznetsov10} requires some further explanation, especially since it
involves an abuse of notation. For a vector with components $X^j$ the quantity
$Y^i=\sum_{j=1}^nB^i_jX^j$ is supposed to be a solution of the equation
$\sum_{j=1}^nA^i_jY^j=X^i$. It is unique up to addition of an element of the
kernel of $A$, i.e. a multiple of $w$. The value of $c$ is independent of
which solution is chosen because of the fact that $a=0$. The matrix with
components $A^i_j$ is not invertible but under the assumptions of the theorem
the vector with components $X^i$ satisfies $\sum_{i=1}^nX^iv_i=0$ and hence lies
in its image.

The bifurcation of Theorem 1 is related to the pitchfork bifurcation. As in the
case of the transcritical bifurcation there is a choice to be made, whether or
not the steady state moves as the parameter $\alpha_1$ is changed. In Theorem 1
it is assumed that it does not move. A corresponding theorem about the
pitchfork bifurcation in one dimension is Theorem 3.17 of \cite{homburg20}. A
version where the steady state is allowed to move is given by Theorem 3.18 of
\cite{homburg20}. The case where the dimension is general is treated in
Theorem 3.21 of \cite{homburg20}. This theme is also discussed in Chapter 1
of \cite{golubitsky85}. Note that the pitchfork bifurcation only involves
one parameter while Theorem 1 involves two. In fact the bifurcation of
Theorem 1 is an unfolding of the pitchfork bifurcation, a concept discussed in
Chapter 3 of \cite{golubitsky85}. An unfolding is a deformation of a bifurcation
depending on additional parameters. Of particular interest are the universal
deformations which include all possible bifurcations close to the original
one. The universal deformation of the pitchfork bifurcation contains three
parameters. There the position of the steady state is allowed to depend on
all parameters. Requiring that it does not depend on the first parameter
leads to a situation similar to that of Theorem 1. Proving that an unfolding
of the pitchfork bifurcation is universal is technically difficult (see
\cite{golubitsky88}, \cite{murdock03}) and trying to adapt the ideas about
universal unfoldings to prove Theorem 1 would, if possible, presumably be
considerably more complicated than the method of proof used above.

We conclude this section with some remarks concerning how this theorem can be
applied in practise. The starting point is a system of equations depending on
many parameters. It is necessary to choose parameters to play the role of
$\alpha_1$ and $\alpha_2$ in the theorem. These parameters should be chosen
so that the coordinates of the DFE do not depend on $\alpha_1$. The quantities
$f_{u\alpha_2}$ and $f_{uu\alpha_1}$ can be rewritten in terms of the full system
as can be done for the quantities $b$ and $d$.

\section{An analysis of Brauer}\label{brauer}

In this section we compare the result of Theorem 1 with an analysis due to
Brauer \cite{brauer04} of a simple vaccination model. Since the results
of \cite{brauer04} are rather complete there is nothing more to be proved about
that model than was done in \cite{brauer04}. The aim of this section is
instead to compare the arguments of \cite{brauer04} with an alternative proof
of those results provided by Theorem 1 of the present paper so as to obtain
insights which can be helpful for the understanding of more complicated models.

The system (4) of Brauer \cite{brauer04} is three-dimensional. By considering
the passage to $\omega$-limit sets he relates it to a two-dimensional system,
(5) of \cite{brauer04}, which he uses to analyse the long-time behaviour of
the original system, in particular the multiplicity and stability of steady
states. The reduced system (5) of Brauer \cite{brauer04} is
\begin{eqnarray}
&&I'=\beta[K-I-(1-\sigma)V]I-(\mu+\gamma)I,\label{basic1}\\
&&V'=\phi[K-I]-\sigma\beta VI-(\mu+\theta+\phi)V.\label{basic2}  
\end{eqnarray}
Here $0\le\sigma\le 1$ and all the other parameters are positive. It is
elementary to show that the system (\ref{basic1})-(\ref{basic2}) has a unique
boundary steady state, the disease-free equilibrium (DFE), which is given by
$\left(0,\frac{\phi K}{\mu+\theta+\phi}\right)$.

The linearization of the system (\ref{basic1})-(\ref{basic2}) at
a general point is
\begin{equation}
\left[
{\begin{array}{cc}
   -2\beta I-\mu-\gamma-\beta(1-\sigma) V+\beta K & -\beta(1-\sigma)I \\
   -(\phi+\sigma\beta V) & -(\mu+\theta+\phi+\sigma\beta I)
\end{array}}
\right].
\end{equation}
At the DFE this reduces to
\begin{equation}
\left[
{\begin{array}{cc}
   -\mu-\gamma-\beta(1-\sigma) V+\beta K & 0 \\
   -(\phi+\sigma\beta V) & -(\mu+\theta+\phi)
\end{array}}
\right].
\end{equation}
Both eigenvalues of this matrix are negative if and only if
\begin{equation}
-\mu-\gamma-\beta(1-\sigma) V+\beta K<0.
\end{equation}
Using the expression for $V$ at the DFE shows that this is equivalent to the
condition
\begin{equation}\label{R0}
\frac{\beta K(\mu+\theta+\sigma\phi)}{(\mu+\gamma)(\mu+\theta+\phi)}<1.
\end{equation}
Given that $\phi$ is the rate of vaccination putting $\phi=0$ in the
expression on the left hand side of (\ref{R0}) gives something which plays the
role of the basic reproductive ratio in the absence of vaccination.

Our aim is to compare this with the approach of van den Driessche and
Watmough to backward bifurcations in \cite{vandendriessche02} and Theorem 1
of the present paper. It is not clear whether the formalism of
\cite{vandendriessche02} can be applied to the reduced system and so we will
apply it to the full system, (4) of \cite{brauer04}, which is
\begin{eqnarray}
  &&S'=\Lambda-\beta SI-(\mu+\phi)S+\gamma I+\theta V,\label{full1}\\
  &&I'=\beta SI+\sigma\beta VI-(\mu+\gamma)I,\label{full2}\\
  &&V'=\phi S-\sigma\beta VI-(\mu+\theta)V.\label{full3}
\end{eqnarray}
The system (\ref{full1})-(\ref{full3}) has the DFE
$\left(\frac{(\mu+\theta)\Lambda}{\mu(\mu+\phi+\theta)},
  0,\frac{\phi\Lambda}{\mu(\mu+\theta+\phi)}\right)$. 
To relate the two- and three-dimensional systems we should take
$K=S+I+V=\frac{\Lambda}{\mu}$. It is important for the analysis to be done
later that the DFE does not depend on $\beta$. The linearization of
(\ref{full1})-(\ref{full3}) at the DFE is
\begin{equation}
\left[
{\begin{array}{ccc}
   -(\mu+\phi) & -\beta S+\gamma & \theta \\
   0 & \beta S+\sigma\beta V-(\mu+\gamma)& 0\\
   \phi & -\sigma\beta V & -(\mu+\theta)
\end{array}}
\right].
\end{equation}
One of the eigenvalues of this matrix is $\beta (S+\sigma V)-(\mu+\gamma)$.
The other two are the eigenvalues of a two by two matrix with negative trace and
positive determinant and so they have negative real parts. Note that at the DFE
\begin{equation}
S+\sigma V=\frac{(\mu+\theta+\sigma\phi)\Lambda}{\mu(\mu+\theta+\phi)}
\end{equation}  
so that the condition for all eigenvalues to have negative real parts is the
same as that which we found for the reduced system.

In the case that there is a zero eigenvalue the right eigenvector of the
linearization corresponding to the eigenvalue zero can be taken to be of the
form $(E,1,-1-E)$, with $E=\frac{\gamma-\beta S-\theta}{\mu+\phi+\theta}$.
A corresponding left eigenvector can be taken to be of the form $(0,1,0)$. We
are interested in this quantity at the bifurcation point. There
we have the relation
\begin{equation}
  \frac{\mu+\theta+\sigma\phi}{\mu+\theta}\beta S
  =\beta (S+\sigma V)=\mu+\gamma.
\end{equation}
It follows that
\begin{equation}
  \gamma-\theta-\beta S
=\frac{-(\mu+\theta)^2+\sigma\phi(\gamma-\theta)}{\mu+\theta+\sigma\phi}.
\end{equation}

We now want to apply the method of \cite{vandendriessche02} to this system,
taking the parameter $\mu$ of the general method to be equal to
$\beta-\beta^*$, where $\beta^*$ is the value of $\beta$ at the
bifurcation point. Contact can be made with that formalism by choosing
$x_1=I$, $x_2=S$, $x_3=V$, with $I$ being the only infected variable. Then the
conditions (A1)-(A5) of \cite{vandendriessche02} are satisfied. On the other
hand we choose the order $(S,I,V)$ of the variables as in
(\ref{full1})-(\ref{full3}) in the calculations
which follow. Since the only non-zero component of the left eigenvector is the
second we only need second derivatives of the right hand side of the evolution
equation for $I$ to calculate the coefficient $a$ of \cite{vandendriessche02}.
Up to a positive factor the coefficient $a$ is equal to
$E-(1+E)\sigma=E(1-\sigma)-\sigma$. This means that provided $\sigma<1$ the
condition $a=0$ is equivalent to $E=\frac{\sigma}{1-\sigma}$. In the case
$\sigma=1$ the quantity $a$ is never zero and the following analysis does not
apply. As mentioned in \cite{brauer04} this can be interpreted as the case in
which the vaccine has no effect. In \cite{brauer04} the bifurcation parameter
used is $\beta$. To calculate the coefficient $b$ we only need first derivatives
of the right hand side of the evolution equation for $I$. It is equal to
$S+\sigma V$, which is positive. The condition for a backward bifurcation is
$a>0$, which is equivalent to
\begin{equation}
  \sigma^2\phi^2+\sigma\phi[-\gamma+\mu+2\theta+\sigma(\gamma+\mu)]
  +(\mu+\theta)^2<0.\label{backward1}
\end{equation}

To compare with \cite{brauer04} the equations for steady states of the reduced
system with $I>0$ will be considered. These are
\begin{eqnarray}
&&\beta[K-I-(1-\sigma)V]=\mu+\gamma,\\
&&\phi[K-I]=[\sigma\beta I+(\mu+\theta+\phi)]V.
\end{eqnarray}
The first of these equations implies that
$\beta(1-\sigma)V=\beta(K-I)-\mu-\gamma$. Hence multiplying the second by
$\beta(1-\sigma)$ gives
\begin{equation}
\beta(1-\sigma)\phi[K-I]=[\sigma\beta I+(\mu+\theta+\phi)]
[\beta(K-I)-\mu-\gamma].
\end{equation}
This can be rewritten in the form $a_2I^2+a_1I+a_0=0$, where
\begin{eqnarray}
  &&a_2=\sigma\beta^2,\\
  &&a_1=\sigma\beta(\mu+\gamma)+\beta(\mu+\theta+\sigma\phi)-\sigma\beta^2K,\\
  &&a_0=(\mu+\theta+\phi)(\mu+\gamma)-\beta(\mu+\theta+\sigma\phi)K.
\end{eqnarray}
Note that $a_0=0$ precisely when the condition for a zero eigenvalue is
satisfied. Then there is a bifurcation of the steady states. It gives rise
to a backward bifurcation in the biologically relevant region precisely
when the minimum of the parabola occurs at a positive value of $I$. The
minimum occurs when $2a_2I+a_1=0$, i.e. when $I=-\frac{a_1}{2a_2}$. Since
$a_2>0$ the condition for a backward bifurcation is $a_1<0$.

Now $K$ will be eliminated from the system of equations $a_0=a_1=0$.
Multiplying the equation $a_0=0$ by $\sigma\beta$, substituting the equation
$a_1=0$ into this and sorting according to powers of $\phi$ gives
\begin{equation}
  \beta\sigma^2\phi^2+\beta\sigma\phi[-\gamma+\mu+2\theta
  +\sigma(\mu+\gamma)]+\beta (\mu+\theta)^2=0.\label{backward2}
\end{equation}
Comparing (\ref{backward1}) and (\ref{backward2}) makes manifest that the
conditions for a backward bifurcation obtained by the two methods are
equivalent. Note that in the exceptional case $\sigma=1$ no backward
bifurcation is possible.

Next we compare the analysis of \cite{brauer04} with Theorem 1. In the context
of that theorem we introduced a coefficient $c$. Now this quantity
will be considered for the system (\ref{full1})-(\ref{full3}). Since the right
hand sides of these equations are all polynomials which are at most quadratic
the third order derivatives make no contribution in this case. Thus we only
need to consider the second order derivatives. A first step towards doing this
calculation is to evaluate the vector with components
$y_l=\sum_{j,k=1}^3\frac{\partial^2 f_l}{\partial x_j\partial x_k}w_jw_k$.
It is of the form $2\beta[-E,-\sigma+E(1-\sigma),\sigma (1+E)]$.
We now want to find a vector whose image under the linearization at the
bifurcation point is the vector with components $y_l$. This is made possible
by the fact that $y_2$ is proportional to $a$ and hence zero in this case.
The vector with components $y_l$ can be rewritten as
$\frac{2\beta}{1-\sigma}[-\sigma,0,\sigma]$. The non-uniqueness of the vector
with this image can be exploited to suppose that it is of
the form $(B,0,C)$ and then we can compute that
\begin{eqnarray}
  &&B=\frac{2\beta\sigma}{(1-\sigma)(\phi+\mu+\theta)},\\
  &&C=-\frac{2\beta\sigma}{(1-\sigma)(\phi+\mu+\theta)}.
\end{eqnarray}
This gives
\begin{equation}
  c=-\frac{\partial^2f_2}{\partial S\partial I}B
  -\frac{\partial^2f_2}{\partial I\partial V}C
  =-\frac{2\beta^2\sigma}{\phi+\mu+\sigma}.
\end{equation}
Thus $c$ is negative for all values of the parameters.
To apply Theorem 1 to this system it is necessary to make a choice of the
parameters $\alpha_1$ and $\alpha_2$ which ensures that $e\ne 0$. This can
be achieved by choosing $\alpha_1=\beta-\beta^*$ and $\alpha_2=\sigma-\sigma^*$,
where $\sigma^*$ is the value of $\sigma$ at the bifurcation point. Then
$f_{uu\alpha_1}(0,0,0)=\beta a=0$, $f_{u\alpha_1}(0,0,0)>0$ and
$f_{uu\alpha_2}(0,0,0)<0$. It follows that in this case $e>0$ and Theorem 1
implies the existence of multiple steady states in cases with ${\cal R}_0<1$
but does not imply the existence of multiple steady states in cases with
${\cal R}_0>1$. In this way it reproduces the conclusions of \cite{brauer04}
by a different method.

The models considered in this section were originally introduced in
\cite{kribs00}. That paper introduced both the three- and two-dimensional
models but did not provide a precise explanation of the relation between them.
This issue is clarified in \cite{brauer04}. The authors of \cite{kribs00} carry
out a centre manifold analysis of the two-dimensional model and this allows
them to establish the stability properties of the steady states. To do this
they make use of the fact that since steady states correspond to roots of a
quadratic equation the cusp bifurcation point can be determined explicitly for
this model. This is used in \cite{kribs00} to do a centre manifold analysis at
the cusp bifurcation point.

\section{A model of Martcheva}\label{martcheva}

In this section we examine the applicability of Theorem 1 to a model system
for backward bifurcations studied by Martcheva \cite{martcheva19}
(Example 2.1 of that paper). It is a description of cholera with
vaccination. That system shares the property of the model of Brauer
considered in the last section that its steady states are in one-to-one
correspondence with the roots of a quadratic equation. An analysis related to
the method of van den Driessche and Watmough
\cite{vandendriessche02} is carried out for
this model in \cite{martcheva19}. (Instead of referring primarily to
\cite{vandendriessche02} it uses the later work of \cite{castillo04} as a
basis.) The formulation with the quadratic equation provides information
about the number of steady states in this model. Other calculations in
\cite{martcheva19} give the type of information about stability which can
be obtained using the method of \cite{vandendriessche02}. On the other hand the
results of \cite{martcheva19} do not suffice to prove that there is a choice of
the parameters in this model for which there exist one stable and one unstable
positive steady state. For the model of Brauer considered in the last section
information of that type was obtained in \cite{brauer04} but this made use of
the fact that the analysis could be reduced to that of a two-dimensional model
in that case.

The model of \cite{martcheva19} is 
\begin{eqnarray}
  &&S'(t)=\Lambda-\beta\frac{S(t)B(t)}{B(t)+D}-(\mu+\psi)S(t)+wR(t),
  \label{martcheva1}\\
  &&V'(t)=\psi S(t)-\sigma\beta\frac{V(t)B(t)}{B(t)+D}-\mu V(t),
  \label{martcheva2}\\
  &&I'(t)=\beta\frac{S(t)B(t)}{B(t)+D}+\sigma\beta\frac{V(t)B(t)}{B(t)+D}
  -(\mu+\gamma)I(t),\label{martcheva3}\\
  &&R'(t)=\gamma I(t)-(\mu+w)R(t),\label{martcheva4}\\
  &&B'(t)=\eta I(t)-\delta B(t).\label{martcheva5}
\end{eqnarray}  
The unique DFE is $\left(\frac{\Lambda}{\mu+\psi},
  \frac{\Lambda\psi}{\mu(\mu+\psi)},0,0,0\right)$. It is important for the
analysis to be done later that the DFE does not depend on $\eta$.
In \cite{martcheva19} the basic reproductive ratio for this model (in the
presence of vaccination) has been given as
${\cal R}_0=\frac{\eta\beta(S+\sigma V)}{D\delta (\mu+\gamma)}$ and
the coefficients $a$ and $b$ which are important for backward bifurcations are
computed there.  The linearization at the disease-free equilibrium and some of
the second order derivatives are also given. With a suitable normalization the
left eigenvector of the linearization with eigenvalue zero is
$\left(0,0,\frac{\eta}{\mu+\lambda},0,1\right)$ and the components $w_3$,
$w_4$ and $w_5$ of the corresponding right eigenvector are all positive. To
relate this system to the formalism of \cite{vandendriessche02} we should
choose $I$ and $B$ to be the infected variables.

One way to show that there exist parameters for which this model has one stable
and one unstable steady state would be to apply Theorem 1 and we now want to
investigate this possibility. Note first that the case of Theorem 1 with $c>0$
cannot occur for this model.

\noindent
{\bf Lemma 1} Consider a bifurcation point in the parameter space of
(\ref{martcheva1})-(\ref{martcheva5}) for which ${\cal R}_0=1$, $a=0$,
$b>0$ and $e\ne 0$. Then $c\le 0$.

\noindent
{\bf Proof} For this system the steady states are in one-to-one correspondence
with the roots of a quadratic equation and it can be shown that the number of
solutions of this equation counting multiplicity is even for ${\cal R}_0<1$
and odd for ${\cal R}_0>1$. Under
the assumptions of Lemma 1 if $c$ is positive Theorem 1 implies that
there exist parameter values with ${\cal R}_0>1$ for which there exist at
least two steady states. Since the steady states are in one-to-one
correspondence with the roots of a quadratic equation there are at most two
of them and hence in this case exactly two. This contradicts the statement
about the parity of the number of roots made above. It follows that
the assumption $c>0$ leads to a contradiction and that in fact $c\le 0$.
$\blacksquare$.

To show that Theorem 1 can be applied we need to know that the coefficient $c$
is negative for some choice of the parameters. First we calculate the components
$y_l=\sum_{j,k=1}^5\frac{\partial^2 f_l}{\partial x_j\partial x_k}w_jw_k$.
Clearly $y_4=y_5=0$. The other components are
\begin{eqnarray}
  &&y_1=\frac{2\beta}{D}\left(-w_1
  +\frac{S}{D}w_5\right)w_5=\frac{2\beta}{D}\left(-w_1
  +\frac{S}{D}\right),\label{y1}\\
  &&y_2=\frac{2\beta\sigma}{D}\left(-w_2+\frac{V}{D}w_5\right)w_5
     =\frac{2\beta\sigma}{D}\left(-w_2+\frac{V}{D}\right),\label{y2}\\
  &&y_3=\frac{2\beta}{D}\left[w_1+\sigma w_2
     -\frac{(S+\sigma V)}{D}w_5\right]w_5.\label{y3}
\end{eqnarray}
It follows from (\ref{y2}) and \cite{martcheva19} that $y_3=2\beta w_5a$.
Hence in the case $a=0$ we have $y_3=0$. When looking for a vector whose image
under the linearization is equal to $y_l$ it is possible and convenient to
choose the last component to be zero. It then follows using the fact that
$y_3=y_4=0$ that the third and fourth components are also zero. For the
remaining components we have

\begin{eqnarray}
  &&\Sigma_{l=1}^5B_l^1y_l=-\frac{1}{\mu+\psi}y_1
     =-\frac{2\beta}{(\mu+\psi)D}\left(-w_1+\frac{S}{D}\right),\\
  &&\Sigma_{l=1}^5B_l^2y_l=-\frac{1}{\mu}\left(\frac{\psi}{\mu+\psi}y_1+y_2\right)
     \nonumber\\
  &&=-\frac{2\beta\sigma}{\mu D}\left[\frac{\psi}{\mu+\psi}
     \left(-w_1+\frac{S}{D}\right)+\sigma\left(-w_2+\frac{V}{D}\right)
     \right].
\end{eqnarray}     
The contribution to $c$ from the second derivatives, call it $c^{(2)}$, is
\begin{eqnarray}
  &&-\Sigma_{i,j,k,l=1}^5B_l^jy_l
     \frac{\partial^2 f_i}{\partial x_j\partial x_k}v_iw_k\\
  &&=-\Sigma_{k,l=1}^5 \left[B_l^1y_l
     \frac{\partial^2 f_3}{\partial x_1\partial x_5}
     w_1+B_l^2y_l\frac{\partial^2 f_3}{\partial x_2\partial x_5}
     w_2\right]v_3.
\end{eqnarray}
The second derivatives in this expression which are non-zero are
$\frac{\partial^2f_3}{\partial S\partial B}=\frac{\beta}{D}$ and
$\frac{\partial^2f_3}{\partial V\partial B}=\frac{\sigma\beta}{D}$.
Substituting this information into the above expression gives
\begin{equation}
c^{(2)}=-\frac{\beta}{D}\left(\sum_{l=1}^5B_l^1y_lw_1
    +\sigma\sum_{l=1}^5 B_l^2y_lw_2\right)v_3.
\end{equation}
More explicitly this is
\begin{equation}
  \frac{2\beta^2}{D^2}\left\{\frac{1}{\mu+\psi}
    \left(-w_1+\frac{S}{D}\right)w_1
    +\frac{\sigma}{\mu}\left[\frac{\psi}{\mu+\psi}\left(-w_1+\frac{S}{D}\right)
      +\sigma\left(-w_2+\frac{V}{D}\right)\right]w_2\right\}
\end{equation}
We want to evaluate this at the DFE.

The only relevant third derivatives are
$\frac{\partial^3f_3}{\partial S\partial B^2}=-\frac{2\beta}{D^2}$,
$\frac{\partial^3f_3}{\partial V\partial B^2}=-\frac{2\sigma\beta}{D^2}$ and
$\frac{\partial^3f_3}{\partial B^3}=\frac{6\beta (S+\sigma V)}{D^3}$
and so the contribution of the third derivatives to $c$ is
\begin{eqnarray}
&&  c^{(3)}=\frac{\beta}{3D^2} v_3\left(-6(w_1+\sigma w_2)+\frac{6(S+\sigma V)}{D}
  \right)\nonumber\\
&&  =\frac{2\beta}{D^2} v_3\left[-(w_1+\sigma w_2)
    +\frac{\Lambda}{D(\mu+\psi)}\left(1+\frac{\psi\sigma}{\mu}\right)
  \right].
\end{eqnarray}
Note that from \cite{martcheva19}
\begin{equation}
  a=2\beta\left(\frac{1}{D}w_1+\sigma\frac{1}{D}w_2
  -\frac{S+\sigma V}{D^2}\right).\label{azero}
\end{equation}
Thus the condition $a=0$ is equivalent to
\begin{equation}\label{w1w2}
w_1+\sigma w_2=\frac{\Lambda}{D(\mu+\psi)}\left(1+\frac{\sigma\psi}{\mu}\right)
\end{equation}
so that in this case $c^{(3)}=0$.

If a positive steady state is given for any values of the parameters then we
can obtain a corresponding steady state with ${\cal R}_0=1$ by rescaling
$\eta$.

\noindent
{\bf Theorem 2} There exist parameter values for the system
(\ref{martcheva1})-(\ref{martcheva5}) for which ${\cal R}_0<1$ and there
exist exactly two positive steady states, one stable and one unstable.

\noindent
{\bf Proof} This will be proved by applying Theorem 1 with a suitable
choice of parameters. We choose $\alpha_1=\eta$ and
$\alpha_2=D$. We first claim that there exists a bifurcation point with
${\cal R}_0=1$, $a=0$, and $b>0$. It was shown in \cite{martcheva19} that
the coefficient $b$ is always positive. Note however that the definition
of $b$ depends on the parameter chosen when applying the method of
\cite{vandendriessche02}. In \cite{martcheva19} the parameter
chosen was $\beta$ while we will instead choose $\eta$. However $b$ is also
always positive for the latter choice. Finding a point satisfying the first
two conditions corresponds to choosing parameters so that the coefficients
$a_0$ and $a_1$ in the quadratic equation in \cite{martcheva19} both vanish.
To find parameters of this kind we start with the expression on the left hand
side of equation (14) in \cite{martcheva19} whose sign is equal to that of
$a_1$ under the condition that $a_0=0$. It is
\begin{equation}
  -\frac{S^0}{D}Q-\frac{\beta^* V^0}{D}\frac{\sigma^2}{\mu}
  -\frac{\beta^*S^0}{D(\mu+\psi)}Q+\frac{q}{\mu+\psi}Q.
\end{equation}
After substituting for the quantities $S^0$, $V^0$, $\beta^*$, $q$ and $Q$ using
the expressions given in \cite{martcheva19} and simplifying this becomes
\begin{equation}\label{a0a1}
  -\frac{\Lambda(\mu+\sigma\psi)}{D\mu(\mu+\psi)}
  -\frac{\delta (\mu+\gamma)}{\eta(\mu+\sigma\psi)}\frac{\psi\sigma^2}{\mu}
  +\frac{\delta}{\eta(\mu+\psi)}\left[-\mu-\gamma
  +\frac{w}{\mu+w}\frac{\gamma (\mu+\sigma\psi)}
  {\mu}\right].
\end{equation}
It will now be shown that this quantity can take on both signs. Making $\Lambda$
sufficiently large while keeping the other parameters fixed makes this quantity
negative. To see that there are parameters which make it positive first choose
$\psi=\frac{2\mu}{\sigma}\frac{\mu+w}{w}$. Then the expression in square
brackets in (\ref{a0a1}) becomes $\gamma\left(1+\frac{w}{\mu+w}\right)-\mu$.
Choosing $\gamma$ large enough makes this expression positive. Making $\sigma$
sufficiently small while fixing the value of $\sigma\psi$ and those of the
other parameters allows the second term in (\ref{a0a1}) to be dominated by that
in the square brackets. Finally, making $\Lambda$ sufficiently small allows
the first term in (\ref{a0a1}) to be dominated by that in square brackets so
that it can be arranged that $a_1>0$.

With the given choice of $\alpha_1$ and $\alpha_2$ we have $e\ne 0$.
For as in the case of the system of Brauer we have $f_{uu\alpha_1}=0$,
$f_{u\alpha_1}>0$ and $f_{uu\alpha_2}\ne 0$. In fact the last quantity is a
non-zero multiple of $w_1+\sigma w_2$ and so the fact that it is non-zero
follows from (\ref{w1w2}). If at the given point $c<0$ then
Theorem 1 can be applied and the proof of Theorem 2 is complete. Thus it
suffices to show that the assumption that $c$ is always non-negative leads to a
contradiction. Suppose that $c\ge 0$ in a open neighbourhood of the given
point in parameter space. Combining this with Lemma 1 implies that $c=0$ in
that neighbourhood. This means that the expression in (\ref{a0a1}) vanishes
in that neighbourhood. This gives a contradiction since, for instance,
changing $\Lambda$ while keeping the other parameters fixed changes the value
of that expression. $\blacksquare$

\section{An in-host model for hepatitis C}\label{hepc}

Consider the following in-host model for hepatitis C introduced in
\cite{nangue22}.
\begin{eqnarray}
  &&\frac{dT}{dt}=s+r_TT\left(1-\frac{T+I}{T_{\rm max}}\right)-dT
     -\frac{bTV}{T+I}\label{basichepc1},\\
  &&\frac{dI}{dt}=r_II\left(1-\frac{T+I}{T_{\rm max}}\right)
  +\frac{bTV}{T+I}-\delta I,\label{basichepc2}\\
  &&\frac{dV}{dt}=\rho R^*I-cV-\frac{bTV}{T+I}.\label{basichepc3}
\end{eqnarray}
It can be related to the class of systems considered in Section \ref{main}
with $n=3$ and $m=2$. In that correspondence $I$ and $V$ are the infected
variables and correspond to $x_1$ and $x_2$ while $T$ corresponds to $x_3$ and
has the value $p_0$ at the DFE $x_0$, where $p_0$ is as in \cite{nangue22}. In
the calculations which follow we order the unknowns in a different way, as in
(\ref{basichepc1})-(\ref{basichepc3}). Thus the ordering is $(T,I,V)$ and the
DFE is $(p_0,0,0)$. It is important for the analysis to be done later that the
DFE does not depend on $\rho$. From this point on the quantities $b$,
$c$ and $d$ introduced in Section \ref{main} are denoted by $b_{VW}$,
$c_{VW}$ and $d_{VW}$ to prevent them being confused with the parameters $b$,
$c$ and $d$ occurring in the system (\ref{basichepc1})-(\ref{basichepc3}).
The vectors $v$ and $w$ in the general analysis of \cite{vandendriessche02} 
were computed for this example in \cite{nangue23}. To give explicit
expressions for these quantities we recall the following notation from
\cite{nangue22}.
\begin{eqnarray}
  &&p_0=\left(r_T-d+\sqrt{(r_T-d)^2+\frac{4sr_T}{T_{\rm max}}}
  \right)\frac{T_{\rm max}}{2r_T},\\
  &&a_{11}=\sqrt{(r_T-d)^2+\frac{4sr_T}{T_{\rm max}}},\\
  &&a_{12}=\frac{p_0r_T}{T_{\rm max}},\\
  &&a_{22}=-\left[\delta+r_I\left(\frac{p_0}{T_{\rm max}}-1\right)\right].     
\end{eqnarray}
Then $v=(0,b+c,b)$ and $w=(-a_{12}(b+c)-b\rho R^*,a_{11}(b+c),a_{11}\rho R^*)$.
Note that $p_0>0$ and $0<a_{12}<a_{11}$. In this model
${\cal R}_0=\frac{b\rho R^*}{(b+c)(\delta-r_I(1-\frac{p_0}{T_{\rm max}}))}$.
From now on we only consider sets of parameters for which
$\delta>r_I(1-\frac{p_0}{T_{\rm max}})$, so that ${\cal R}_0>0$. In that case
$a_{22}<0$. Suppose that $p_0\le T_{\max}$. Then
given any parameter set satisfying ${\cal R}_0>0$ it is possible to
modify $\delta$ while keeping all other parameters fixed so that
${\cal R}_0=1$. More specifically, we choose
\begin{equation}\label{fixdelta}
  \delta=\frac{b\rho R^*}{b+c}
    +r_I\left(1-\frac{p_0}{T_{\rm max}}\right).
\end{equation}  
We refer to the resulting parameter set as the corresponding
parameter set with ${\cal R}_0=1$.

In order to check the criteria of Theorem 4 of \cite{vandendriessche02} and 
Theorem  1 of Section 2 we first need to calculate the derivatives of the
right hand side of the equations at the bifurcation point up to order three. The
linearization at the bifurcation point is (cf. \cite{nangue22})
\begin{equation}
\left[
{\begin{array}{ccc}
   -a_{11} & -a_{12} & -b \\
   0 & a_{22} & b \\
   0 & \rho R^* & -c-b
\end{array}}
\right].
\end{equation}
Consider next the second order derivatives. Most of the expressions for these
derivatives given in \cite{nangue23} are correct but there are two errors
there. The expression for $\frac{\partial^2f_2}{\partial I\partial V}$ should
have a factor $T$ in the numerator  instead of the factor $I$. In addition the
factor $\beta$ in that expression should be $b$. When we evaluate at the
disease-free steady state we see that the only second derivatives which make
non-trivial contributions to $a$ are
$\frac{\partial^2 f_2}{\partial T\partial I}$,
$\frac{\partial^2 f_2}{\partial I^2}$,
$\frac{\partial^2 f_2}{\partial I\partial V}$ and
$\frac{\partial^2 f_3}{\partial I\partial V}$. They are equal to
$-\frac{r_I}{T_{\rm max}}$,
$-\frac{2r_I}{T_{\rm max}}$, $-\frac{b}{p_0}$ and
$\frac{b}{p_0}$, respectively. Here we have used the fact that the derivatives
of $f_1$ make no contribution to $a$ since $v_1=0$.

In later calculations we
will require the second order derivatives of $f_1$ at the bifurcation point
and so we record them here. Those terms which are of degree zero or one
make no contribution and so we only need to consider the other two terms.
Derivatives involving $V$ of the first of these terms are zero. On the other
hand the only derivatives involving $V$ of the second term which make a non-zero
contribution are those containing exactly one derivative with respect to $V$.
It follows that the only relevant derivatives are
$\frac{\partial^2f_1}{\partial T^2}$,
$\frac{\partial^2f_1}{\partial T\partial I}$ and
$\frac{\partial^2f_1}{\partial I\partial V}$. They are equal to
$-\frac{2r_T}{T_{\rm max}}$, $-\frac{r_T}{T_{\rm max}}$ and $\frac{b}{p_0}$,
respectively.

It follows from the expressions for the second derivatives that 
\begin{eqnarray}
&&a=w_2\left[-v_2(w_1+w_2)\frac{r_I}{T_{\rm max}}
   +(-v_2+v_3)w_3\frac{b}{p_0}\right]\nonumber\\
&&=w_2\left[-(b+c)(w_1+w_2)\frac{r_I}{T_{\rm max}}
   -cw_3\frac{b}{p_0}\right]\nonumber\\
&&=a_{11}(b+c)\left[-(b+c)(w_1+w_2)\frac{r_I}{T_{\rm max}}
-cw_3\frac{b}{p_0}\right]\nonumber\\
&&=a_{11}(b+c)\left[(b+c)((b+c)(-a_{11}+a_{12})+b\rho R^*)\frac{r_I}{T_{\rm max}}
-b\rho R^*\frac{ca_{11}}{p_0}\right].\nonumber
\end{eqnarray}
This corrects the expression for $a$ given in \cite{nangue23}. $p_0$, $a_{11}$ 
and $a_{12}$ only depend on $d$, $r_T$, $s$ and $T_{\rm max}$. Note that the
quantity ${\cal R}_0$ depends on all the parameters except $d$. In the treatment
of backward bifurcations in \cite{nangue23} the steady state is fixed while
varying ${\cal R}_0$. This can be done by fixing $d$, $r_T$, $s$ and
$T_{\rm max}$ and varying some other parameter. We can for instance vary $\rho$
which occurs as a factor in ${\cal R}_0$ so as to cause ${\cal R}_0$ to
pass through one with non-zero velocity. The following result corrects
a claim made in \cite{nangue23}.

\noindent
{\bf Theorem 3} There exist parameter values for the system
(\ref{basichepc1})-(\ref{basichepc3}) for which the hypotheses of Theorem 4 of
\cite{vandendriessche02} are satisfied with $a>0$ and parameter values for
which they are satisfied with $a<0$. Thus this system exhibits both backward
and forward bifurcations.

\noindent
{\bf Proof} Recall that $a_{11}$ and $p_0$ are both positive. Let
$F=\frac{(b+c)r_Ip_0}{a_{11}T_{\rm max}}-c$. If $F>0$ then the net contribution
to $a$ of terms containing $\rho$ is positive. Then making $\rho$ large
enough while fixing all other parameters ensures that $a$ is
positive. This quantity should be evaluated at a point where ${\cal R}_0=1$.
We start with a set of parameters for which $F>0$. By choosing $\delta$
according to (\ref{fixdelta}) we can ensure that ${\cal R}_0=1$. This change
does not affect $F$. Next we want to increase $\rho$ while keeping
${\cal R}_0=1$. This can be done by compensating the increase in $\rho$ by a
change in $\delta$.
%All this is consistent with the
%condition $r_I<r_T$, which we might like to have for biological reasons.
To show that there is a generic backward bifurcation we must check that the
coefficient $b_{VW}$ is non-zero. It is equal to $bR^*a_ {11}(b+c)>0$. The
expression for $b_{VW}$ in \cite{nangue23} is not correct but the error makes
no difference to the statement of positivity. In this way we have proved that
backward bifurcations do occur in this model if there exist parameter values
for which $F>0$. If $F<0$ then the net contribution of terms containing $\rho$
to $a$ is negative. In this case we can proceed just as in the case $F>0$ to
get a forward bifurcation.

Thus to prove the theorem we only need to know that the sign of $F$ can be
prescribed arbitrarily by a suitable choice of parameters. To make $F$ positive
it suffices to make $c$ small while fixing the other parameters. This completes
the proof that backward bifurcations occur in the model
(\ref{basichepc1})-(\ref{basichepc3}). To see that $F$ can be made negative
first rewrite it in the form
$F=\frac{br_Ip_0}{a_{11}T_{\rm  max}}
+c\left(\frac{r_Ip_0}{a_{11}T_{\rm  max}}-1\right)$. By making $r_I$ sufficiently
small we can ensure that the second summand is negative. Then making $c$
large enough ensures that $F$ is negative. Thus there is a choice of
parameters for which a forward bifurcation occurs. $\blacksquare$

%In fact we can exhibit an explicit set of parameters for which a backward
%bifurcation occurs. To do this set all the parameters in
%(\ref{basichepc1})-(\ref{basichepc3}) except $\delta$ equal to one and
%$\delta=\frac12$. Then $p_0=1$, $a_{11}=2$ and $a_{12}=1$. It follows that
%${\cal R}_0=1$, so that the DFE is a bifurcation point and that $a=-8$, so
%that a backward bifurcation is obtained.

Backward bifurcations were found by simulation in a model closely related to
(\ref{basichepc1})-(\ref{basichepc3}) in \cite{debroy11}. Apart from differences
of notation that model, (2.1) of \cite{debroy11} differs from that considered
here in two ways. The first is that in \cite{debroy11} mass action kinetics
are used instead of the standard incidence function. The other is that in
\cite{debroy11} the effect of absorption of virions, expressed by the terms
containing $b$ in (\ref{basichepc3}) is omitted. In \cite{debroy11} the authors
showed analytically that there was a necessary condition for a backward
bifurcation in their model which when expressed in our notation is $r_I>r_T$.
It plays no role in our analysis. This means that our model gives a
significantly different biological conclusion. Backward bifurcations can occur
in the case that the proliferation rate of uninfected cells is greater than that
of infected cells.

Now consider the derivatives of order three. Most expressions on the right
hand side of the evolution equations are polynomial of degree at most two
and so make no contribution to the third order derivatives. The only summands
which do make a contribution are those proportional to $\frac{bTV}{T+I}$.
If a third order derivative does not contain a derivative with respect to
$V$ then it vanishes at the bifurcation point. The same is true of a third order
derivative which contains at least two derivatives with respect to $V$.
Thus the only ones which can make a non-zero contribution are those with
one derivative with respect to $V$ and two with respect to the other two
variables. It follows that the derivatives of interest are
$\frac{\partial^3}{\partial I^2\partial V}\left(\frac{bTV}{T+I}\right)
=\frac{2bT}{(T+I)^3}$,
$\frac{\partial^3}{\partial I\partial T\partial V}\left(\frac{bTV}{T+I}\right)
=\frac{b(T-I)}{(T+I)^3}$
and $\frac{\partial^3}{\partial T^2\partial V}\left(\frac{bTV}{T+I}\right)
=\frac{-2bI}{(T+I)^3}$. Evaluating these expressions at the bifurcation point
gives $\frac{2b}{p_0^2}$, $\frac{b}{p_0^2}$ and zero.

The equation $a=0$ can be written in the form
\begin{equation}\label{rI}
  r_I=\frac{T_{\rm  max}}{p_0}\frac{b\rho R^*ca_{11}}{(b+c)[(b+c)(-a_{11}+a_{12})+b\rho R^*]}.
\end{equation}
Note that for some values of the parameters the denominator in this expression
is negative so that $a=0$ is not possible in that case.

Next we apply Theorem 1 to this model with $\alpha_1=\rho-\rho^*$. Here $\rho^*$
is the value of $\rho$ at the bifurcation point. Suppose that
with the given choice of parameters we have
a bifurcation point with $a=0$. The condition $b_{VW}>0$ is satisfied. This
can be thought of as a degenerate backward bifurcation. In order that it not
be too degenerate we need to ensure that $c_{VW}\ne 0$. The contribution to
$c_{VW}$ from the third order derivatives is
\begin{equation}
c_{VW}^{(3)}=\frac2{p_0^2}bc(b+c)a_{11}^2\rho R^*[(a_{11}-a_{12})(b+c)-b\rho R^*].
\end{equation}
The contribution of the second order derivatives is
$-\sum_{j,k,l,q=1}^3B_l^q\frac{\partial^2 f_i
}{\partial x_p\partial x_q}v_iw_py_l$
where $y_l=\sum_{i,j=1}^3\frac{\partial^2 f_l}{\partial x_j\partial x_k}w_jw_k$.
In more detail
\begin{eqnarray}
  &&y_1=-\frac{2r_T}{T_{\rm max}}w_1^2-\frac{2r_T}{T_{\rm max}}w_1w_2
     +\frac{2b}{p_0}w_2w_3,\\
  &&y_2=-\frac{2r_I}{T_{\rm max}}w_2^2-\frac{2r_I}{T_{\rm max}}w_1w_2
     -\frac{2b}{p_0}w_2w_3,\\
  &&y_3=\frac{2b}{p_0}w_2w_3.
\end{eqnarray}
The vector $y$ is in the image of the linearization because it is orthogonal
to the left eigenvector $v$ with eigenvalue zero. The latter statement is not
obvious but a calculation shows that it follows from the condition $a=0$.
Since $w_2\ne 0$ there exists a vector of the form $(B,0,C)$ whose image under
the linearization is the vector with components $y_l$. It follows that
\begin{eqnarray}
  &&B=-a_{11}^{-1}\left(y_1-\frac{b}{b+c}y_3\right)\nonumber\\
  &&   =a_{11}^{-1}\left[\frac{2r_T}{T_{\rm max}}w_1^2
+\frac{2r_T}{T_{\rm max}}w_1w_2-\frac{2bc}{(b+c)p_0}w_2w_3\right],\label{B}\\
  &&C=-\frac{1}{b+c}y_3=-\frac{2b}{(b+c)p_0}w_2w_3.\label{C}
\end{eqnarray}     
The contribution to $c_{VW}$ from the second derivatives is
\begin{eqnarray}
  c_{VW}^{(2)}=(b+c)\left[(b+c)\frac{r_I}{T_{\rm max}}a_{11}B
  +\frac{bca_{11}}{p_0}C\right].
  \label{c2vw}
\end{eqnarray}
Substituting (\ref{B}) and (\ref{C}) into this gives the following expression
for the expression in square brackets in (\ref{c2vw}).
\begin{eqnarray}
&&(b+c)\frac{2r_Ir_T}{T_{\rm max}^2}w_1(w_1+w_2)
   -\frac{2bc}{p_0}
   \left[\frac{r_I}{T_{\rm max}}+\frac{ba_{11}}{(b+c)p_0}\right]w_2w_3
   \label{c2vwexpl}\\
  &&=(b+c)\frac{2r_Ir_T}{T_{\rm max}^2}
     [-a_{12}(b+c)-b\rho R^*][(a_{11}-a_{12})(b+c)-b\rho R^*]\nonumber\\
  &&-\frac{2bc}{p_0}
     \left[\frac{r_I}{T_{\rm max}}+\frac{ba_{11}}{(b+c)p_0}\right]
     a_{11}^2(b+c)\rho R^*\nonumber   
\end{eqnarray}
When $a=0$ this reduces to
\begin{eqnarray}
&&\frac{2r_T}{T_{\rm max}}
   [a_{12}(b+c)+b\rho R^*]\frac{bc}{p_0}\rho R^*a_{11}
   \nonumber\\
  &&-\frac{2bc}{p_0}
     \left[\frac{r_I}{T_{\rm max}}+\frac{ba_{11}}{(b+c)p_0}\right]
     a_{11}^2(b+c)\rho R^*\label{c2vwa0}
\end{eqnarray}

Since this expression does not depend on $\delta$, if the parameters other than
$\delta$ are varied $\delta$ can be adjusted to ensure that ${\cal R}_0=1$.

\noindent
{\bf Theorem 4} Consider the system (\ref{basichepc1})-(\ref{basichepc3}).
There exist parameter values for which ${\cal R}_0<1$ and there exist
exactly two positive steady states, one stable and one unstable. There is a
variation of parameters under which these two steady states coalesce in a
generic fold bifurcation.

\noindent
{\bf Proof} To obtain this result we apply Theorem 1, with the change of
ordering of the unknowns mentioned at the beginning of this section. Thus
we number the unknowns in the order $(T,I,V)$. For each choice of parameters
there exists a unique DFE. There are corresponding quantities ${\cal R}_0$,
$a$, $c_{VW}$ and $e$. To prove this theorem we show that
there exists a choice of parameters for which ${\cal R}_0=1$, $a=0$, $c_{VW}<0$
and $e\ne 0$. When this has been done Theorem 1 implies that there exist
parameters for which ${\cal R}_0<1$ and there exist at least two positive
steady states, one of which is stable and one unstable. Moreover these steady
states are hyperbolic. It was shown in \cite{nangue23} that positive steady
states are in one-to-one correspondence with roots of a cubic polynomial
$p_3$. Theorem 8 of \cite{nangue23} implies that in the case ${\cal R}_0<1$
the number of roots of the polynomial counting multiplicity is even. (This
theorem requires the assumption that $\rho R^*+r_I-\delta>0$ but Lemma 3 of
\cite{nangue23} implies that if any positive steady state exists this
inequality is satisfied.) It is also the case that the number of roots is at
most three. Hence in this case there are exactly two positive steady
states. In fact $a$, $c_{VW}$ and $e$ are independent of $\delta$. Thus it
suffices to show that the  conditions on $a$, $c_{VW}$ and $e$ can be
satisfied since then the first condition can be obtained by choosing $\delta$
to satisfy (\ref{fixdelta}), as previously explained.

It follows from the above discussion that to complete the proof of 
the theorem it only remains to show that the parameters can be chosen in
such a way that $a=0$, $c_{VW}<0$, $e\ne 0$. We make the parameter choices
$\alpha_1=\rho-\rho^*$ and $\alpha_2=r_I-r_I^*$ where the star denotes the
value at the bifurcation point. We consider the system where the parameters in
(\ref{basichepc1})-(\ref{basichepc3}) are such that ${\cal R}_0=1$ for
$\rho=\rho^*$ and that ${\cal R}_0$ increases when $\rho$ increases. The
bifurcation point to be considered is $(p_0,0,0,\rho^*,r_I^*)$.

Suppose we choose $\rho$ at the bifurcation point while using (\ref{rI})
to determine $r_I$. Of course to allow the existence of a positive steady
state it is necessary to choose the parameters so that the right hand side
of (\ref{rI}) is positive. $c_{VW}^{(3)}$ is always negative. Now consider what
happens when $\rho$ is decreased so that the denominator in (\ref{rI})
approaches zero. Then $r_I\to\infty$. In this limit $c_{VW}^{(2)}$ becomes
negative since the first summand in (\ref{c2vwexpl}) is bounded while the
second tends to $-\infty$. Hence a parameter region is reached where
$c_{VW}<0$. At the bifurcation point $f_{uu\alpha_1}=0$ and $d_{VW}\ne 0$. Thus
$e\ne 0$ and all parts of the theorem except the last are
proved. That the last part holds can be seen by examining the proof of Theorem
1. $\blacksquare$

If it were possible to show that there exists a choice of parameters for which
${\cal R}_0=1$, $a=0$, $c_{VW}>0$ and $e\ne 0$ then Theorem 1 would imply that
there exist parameters for which ${\cal R}_0>1$ and there exist at least two
steady states, one of which is stable and one unstable. Unfortunately we were
not able to decide whether parameters satisfying these conditions do exist.
If it were possible then, as in the proof of Theorem 4, Theorem 8 of
\cite{nangue23} could be used to obtain information about the number of steady
states. In the case ${\cal R}_0>1$ the number of roots of the polynomial $p_3$
counting multiplicity is odd and not greater than three. Hence in fact if there
are parameters for which there are at least two steady states there must be
exactly three.

While looking for parameters of the type just discussed our attention was
drawn to the system obtained from (\ref{basichepc1})-(\ref{basichepc3}) by
setting the parameters $s$ and $d$ equal to zero. We call this the truncated
system. It turns out that the steady states of the truncated system have some
interesting and unusual properties, which will be now be discussed. The
calculations done above also apply to the truncated system. In that case
$p_0=T_{\rm max}$, $a_{11}=a_{12}=r_T$ and the expressions for $c^{(2)}_ {VW}$ and
$c^{(3)}_ {VW}$ simplify considerably. In addition the conditions ${\cal R}_0=1$
and $a=0$ simplify to $\delta=\frac{b\rho R^*}{b+c}$ and
$r_I=\frac{cr_T}{b+c}$. When $a=0$ we have
\begin{equation}
  c^{(2)}_ {VW}=-c^{(3)}_ {VW}
  =\frac{2r_T^2b^2c(b+c)\rho R^*(r_T+\rho R^*)}{T_{\rm max}^2}
\end{equation}
Hence for this system, when $a=0$ it follows that $c_{VW}=0$. Compared to
the full system the truncated system has an additional non-negative steady
state for $r_I>\delta$. It is given by
$\left(0,\frac{T_{\rm max}(r_I-\delta)}{r_I},
  \frac{T_{\rm max}\rho R^*(r_I-\delta)}{cr_I}\right)$.
When $s=0$ the polynomial is of the form $p_3(X)=Xp_2(X)$ for a quadratic
polynomial $p_2$.

When ${\cal R}_0=1$ and $a=0$ the truncated model admits a continuum of
steady states. This corresponds to the case that the polynomial $p_2$
vanishes identically. In terms of the parameter $X$ the steady states are
given by
\begin{eqnarray}
  &&T=\frac{T_{\rm max}X[b(\mu+r_T)X+cr_T-b\mu]}{r_T(bX+c)},\\
  &&I=\frac{T_{\rm max}(1-X)[b(\mu+r_T)X+cr_T-b\mu]}{r_T(bX+c)},\\
  &&V=\frac{\mu T_{\rm max}(1-X)[b(\mu+r_T)X+cr_T-b\mu])}{r_T(bX+c)^2}.
\end{eqnarray}
Note that a continuum of steady states was also found for a special choice of
parameters in a model studied in \cite{hews21} which is closely related to
(\ref{basichepc1})-(\ref{basichepc3}). The two models differ (apart from the
choice of notation) only by the fact that the last term in (\ref{basichepc3})
is absent from the model of \cite{hews21}. If one of the two conditions leading
to the existence of a continuum of steady states is not satisfied there is at
most one positive steady state and is satisfies
\begin{equation}
X\left[b(\delta-\mu)+\frac{r_Ib\mu}{r_T}\right]=\frac{r_Ib\mu}{r_T}-c\delta.
\end{equation}
Let us denote this schematically by $AX=B$. If $A$ and $B$ are non-zero and
have the same sign then there is a unique positive solution for $X$. If they
are non-zero and have opposite signs then there is no such solution. If one of
them is zero and the other non-zero then there is no positive solution for
$X$. If both are zero then the equation is identically satisfied and that is
the case where a continuum of steady states exists. In general when $X$ is
given $T$ can be computed using equation (3.5) of \cite{nangue22}, giving
\begin{equation}
T=\frac{T_{\rm max}\left[b\mu X+(c+bX)(r_I-\delta)\right]}{r_I(c+bX)}.
\end{equation}
Then $I=\frac{T(1-X)}{X}$ and $V=\frac{\mu}{c+bX}I$.

\section{Conclusions and outlook}

In this paper we present a new method for establishing rigorous results on the
number and stability of steady states in epidemiological models. This is based
on a theorem, Theorem 1 of this paper, on the existence of a certain type of
bifurcation, which is an unfolding of the pitchfork bifurcation. Depending on
the sign of a certain coefficient $c_{VW}$ this shows the existence of two
positive steady states, one stable and one unstable when the basic
reproductive ratio ${\cal R}_0$ is less than one in the case $c_{VW}<0$ or
greater than one in the case $c_{VW}>0$. In the former case this corresponds
to pictures obtained by simulations is various cases of backward bifurcations.
In the latter case other methods may be used in certain models to establish
that there are a total of three positive steady states, two stable and one
unstable. In particular, a result about the parity of the number of steady
states plays an important role in the proof. Up to now we have not found an
epidemiological model where the latter case occurs. The former case is
illustrated by the example of an in-host model for hepatitis C which had
previously been studied by the authors in \cite{nangue22} and \cite{nangue23}.
Whether the latter case occurs for that model remains to be determined. It
also remains open whether this method can be used to prove the existence of
more than one positive stable steady state in an epidemiological model. The
only case we have found in the literature where bistability has been proved
for a model of this kind is in \cite{wang06} and that model does not seem to
admit the type of bifurcation whose existence is proved by Theorem 1. From a
technical point of view it has the unusual property that the right hand sides
of the equations are not $C^1$ although they are Lipschitz continuous.

In the model for vaccinations from  \cite{brauer04} considered in Section 2
backward bifurcations and moving fold bifurcations occur, the latter only in
the case $c_{VW}<0$. The model is simple enough that the stability of the
solutions can be established rigorously. In the vaccination model from
\cite{martcheva19} considered in Section 3 backward bifurcations occur but it
is more difficult to prove results about the stability of the steady states.
Theorem 1 of the present paper can be applied to this example and gives
results on stability.

What are the features of models of in-host models for viral infections which
lead to backward bifurcations and the type of bifurcation whose existence is
proved by Theorem 1? Both features occur in the model of \cite{nangue22}. That
model uses standard incidence to model the process of infection. There is a
constant source $s$ of healthy cells. Infected cells proliferate with maximal
rate $r_I$. Healthy cells die at rate $d$. The absorption of virions during
infection is included explicitly in the model. In the proof of Theorem 4 the
fact was used that $a_{11}-a_{12}>0$. If we set $s=d=0$ then this condition
fails to be satisfied and we have already seen that a bifurcation of the type
whose existence is proved by Theorem 1 is not possible in that case. If $s>0$
then the proof does work even if $d=0$. If $s=0$ then it only works under the
additional assumption that $r_T<d$. Thus to obtain this phenomenon in this way
we need either a source of healthy cells or that the death rate of healthy
cells is greater than their maximal proliferation rate. Even when $s=d=0$ both
forward and backward bifurcations are obtained as a consequence of Theorem 3.
If $r_I$ is set to zero then the type of bifurcation whose existence is proved
by Theorem 1 is not possible. On the other hand Theorem 3 and its proof extend
directly to the case $s=d=r_I=0$.

A challenge for the future is to find out whether it is possible to use Theorem
1 to prove that there are parameters for which the system
(\ref{basichepc1})-(\ref{basichepc3}) admits three positive steady states. If
this could be done it would remain to prove by some other method that it can
happen that two of these steady states are stable, so that there is
bistability. It would also be interesting to find out systems of biological
interest which allow an application of Theorem 1 with $c_{VW}>0$. Do there
exist alternative criteria for the presence of a moving fold bifurcation,
which require less extensive calculations than those of Theorem 1?

\vskip 20pt

\noindent{\it Acknowledgement} We are grateful to Ale Jan Homburg for pointing
out to us the relevance of the literature on unfoldings to our results.

\end{document}